\documentclass[a4paper,12pt]{amsart}
\usepackage{graphics}
\usepackage{amssymb}
\usepackage{amsmath}
\usepackage{latexsym}
\usepackage{amscd}
\newtheorem{Th}{Theorem}[section]
\newtheorem{Prop}[Th]{Proposition}
\newtheorem{Def}[Th]{Definition}
\newtheorem{Rem}[Th]{\sc Remark}

\newtheorem{Lemma}[Th]{Lemma}
\newcommand{\be}{\begin{eqnarray*}}
\newcommand{\ee}{\end{eqnarray*}}

\newcommand{\pkXY}{{\mathcal P}(^k\!X,Y)}
\newcommand{\LkXY}{{\mathcal L}(^k\!X,Y)}
\newcommand{\pkucXY}{{\mathcal P}_{\mathop{\rm uc}\nolimits}(^k\!X,Y)}
\newcommand{\pkucmXY}{{\mathcal P}_{\mathop{\rm uc+}\nolimits}(^k\!X,Y)}
\newcommand{\pkuc}{{\mathcal P}_{\mathop{\rm uc}\nolimits}(^k\!}
\newcommand{\pkucm}{{\mathcal P}_{\mathop{\rm uc+}\nolimits}(^k\!}
\newcommand{\pkKXY}{{\mathcal P}_{\mathop{\rm co}\nolimits}(^k\!X,Y)}
\newcommand{\pkK}{{\mathcal P}_{\mathop{\rm co}\nolimits}(^k\!}

\newcommand{\pk}{{\mathcal P}(^k\!}

\newcommand{\LXY}{{\mathcal L}(X,Y)}
\newcommand{\LucXY}{{\mathcal L}_{\mathop{\rm uc}\nolimits}(X,Y)}

\newcommand{\Lk}{{\mathcal L}(^k\!}
\newcommand{\LkX}{{\mathcal L}(^k\!X)}

\newcommand{\Proof}{\noindent {\it Proof. }}
\newcommand{\fin}{\hspace*{\fill} $\Box$}
\newcommand{\finesp}{\hspace*{\fill} $\Box$\vspace{.5\baselineskip}}

\newcommand{\N}{\ensuremath{\mathbb{N}}}

\newcommand{\ra}{\rightarrow}
\newcommand{\lra}{\longrightarrow}
\newcommand{\Ra}{\Rightarrow}       
\newcommand{\ptpp}{\widehat{\otimes}_\pi}
\newcommand{\espv}{\vspace{.5\baselineskip}}

\newcommand{\eps}{\epsilon}
\newcommand{\kdots}{\stackrel{(k)}{\ldots}}
\def\dist{\mathop{\rm dist}\nolimits}
\def\sign{\mathop{\rm sign}\nolimits}

\begin{document}
\title{Orlicz-Pettis polynomials on Banach spaces}

\author[M. Gonz\'alez]{Manuel Gonz\'alez}
\address{Departamento de Matem\'aticas \\
      Facultad de Ciencias\\
      Universidad de Cantabria \\ 39071 Santander (Spain)}
\email{gonzalem@ccaix3.unican.es}
\thanks{The first named author was supported
in part by DGICYT Grant PB 97--0349 (Spain)}

\author[J. M. Guti\'errez]{Joaqu\'\i n M. Guti\'errez}
\address{Departamento de Matem\'atica Aplicada\\
      ETS de Ingenieros Industriales \\
      Universidad Polit\'ecnica de Madrid\\
      C. Jos\'e Guti\'errez Abascal 2 \\
      28006 Madrid (Spain)}
\email{jgutierrez@etsii.upm.es}
\thanks{The second named author was supported
in part by DGICYT Grant PB 96--0607 (Spain)}
\thanks{\hspace*{\fill}\scriptsize file oppbs.tex}



\begin{abstract}
We introduce the class of Orlicz-Pettis polynomials between Banach
spaces, defined by their action on weakly unconditionally Cauchy
series. We give a number of equivalent definitions, examples and
counterexamples which highlight the differences between these
polynomials and the corresponding linear operators.
\end{abstract}

\maketitle

\section{Introduction}

In the study of the isomorphic properties of Banach spaces, some
classes of (bounded linear) operators have been introduced which
include the isomorphisms and preserve certain properties of the
spaces. These are the {\it semigroups of operators,} such as the
semi-Fredholm operators, associated to the ideal of compact
operators
\cite{GO1,GO2}, the tauberian operators, associated to the weakly
compact operators
\cite{KW,GO3}, the {\it Orlicz-Pettis operators,} related to the
unconditionally converging operators \cite{GM}, etc. The
semigroups and the operator ideals are somehow opposite notions:
for every Banach space $X$, the identity map $I_X$ belongs to all
the semigroups, while the null operator belongs to all the ideals
\cite{AGM}.

These semigroups do not have an exact analogue within the class of
polynomials between Banach spaces. Nevertheless, in the present
paper we introduce the class of Orlicz-Pettis polynomials, related
to the unconditionally converging polynomials, and show that they
share certain properties with the Orlicz-Pettis (linear) operators
(see Section~\ref{results}), and do not satisfy some others
(Section~\ref{counterex}). In order to obtain the results of
Section~\ref{counterex}, we are led to give a number of
counterexamples, mainly of vector valued polynomials on $c_0$,
which is the key space when we deal with weakly unconditionally
Cauchy series. These counterexamples are of independent interest
and can give new insight into the differences between linear
operators and polynomials.\espv

Throughout the paper, $X$, $Y$ and $Z$ denote Banach spaces, $X^*$
is the dual of $X$, $B_X$ is its closed unit ball, $\LXY$ stands
for the space of operators  from $X$ into $Y$, $\pkXY$ represents
the space of all $k$-homogeneous (continuous) polynomials from $X$
into $Y$, $\LkXY$ is the space of all $k$-linear (continuous)
mappings from $X^k$ into $Y$. When the range space $Y$ is omitted,
it is supposed to be the scalar field (real or complex). We denote
by $X\ptpp Y$ the projective tensor product of $X$ and $Y$; the
product of $k$ spaces is represented by $\ptpp^k
X:=X\ptpp\cdots\ptpp X$. We use the notation
$x^{(k)}:=x\otimes\kdots\otimes x$, where $x\in X$. The set of
natural numbers is denoted by $\N$, and $(e_n)$ is the unit vector
basis of the space $c_0$. The coordinates of a vector $x\in c_0$
are denoted by $x(i)$ $(i=1,2,\ldots)$.

A formal series $\sum x_n$ in $X$ is {\it weakly unconditionally
Cauchy\/} ({\it w.u.C.,} for short) if, for every $\phi\in X^*$,
we have $\sum|\phi (x_n)|<+\infty$. Equivalent definitions may be
seen in \cite[Theorem~V.6]{Di}. The series is {\it unconditionally
convergent\/} ({\it u.c.,} for short) if every subseries
converges. Equivalent definitions may be seen in
\cite[Theorem~1.9]{DJT}.

It may be helpful to recall that every polynomial between Banach
spaces takes w.u.C. (resp.\ u.c.) series into w.u.C. (resp.\ u.c.)
series \cite[Theorem~2]{GG}. The following simple fact will also
be useful:

\begin{Prop}
Given a polynomial $P\in\pkXY$ and a w.u.C. series $\sum x_n$ in
$X$, the sequence $\left( P\left(\sum_{i=1}^nx_i\right)\right)_n$
is weak Cauchy.
\end{Prop}

\Proof
Let $j\in{\mathcal L}(c_0,X)$ be the operator given by
$j(e_n)=x_n$. Then $P\circ j\in\pk c_0,Y)$. Since $c_0$ has the
Dunford-Pettis property, $P\circ j$ takes weak Cauchy sequences
into weak Cauchy sequences \cite{Ry}. So, the sequence
$$
\left( P\left(\sum_{i=1}^nx_i\right)\right)_n=
\left( P\circ j\left(\sum_{i=1}^ne_i\right)\right)_n
$$
is weak Cauchy.\finesp

A polynomial $P\in\pkXY$ $(k\geq 1)$ is {\it unconditionally
converging\/} \cite{FU,FUExt} if, for each w.u.C. series $\sum
x_n$ in $X$, the sequence $\left( P\left(\sum_{i=1}^n
x_i\right)\right)_n$ is convergent in $Y$. The space of all
unconditionally converging polynomials is denoted by $\pkucXY$ (or
$\LucXY$ if $k=1$). This class of polynomials has been very useful
for obtaining polynomial characterizations of Banach space
properties (see \cite{GV}). Easily, $T\in\LucXY$ if and only if
for each w.u.C. series $\sum x_n$ in $X$, the series $\sum T(x_n)$
is u.c.\ in $Y$. The polynomial $P\in\pkXY$ is {\it (weakly)
compact\/} if $P(B_X)$ is relatively (weakly) compact in $Y$. The
space of compact polynomials from $X$ into $Y$ is denoted by
$\pkKXY$. Every weakly compact polynomial is unconditionally
converging, and every unconditionally converging polynomial on
$c_0$ is compact (see \cite{FU} or \cite{GV}).

The standard notations and definitions in Banach space theory may
be seen in \cite{Di}. For the basics in the theory of polynomials,
we refer to \cite{Din,Mu}.

\section{Positive results}
\label{results}

In this Section, we introduce the Orlicz-Pettis polynomials, as
those satisfying the following main result. We give some other
properties, and a first example of a polynomial in this class.

\begin{Th}
\label{uc+}
Given $k\in\N$ and $P\in\pkXY$, the following assertions are
equivalent:

{\rm (A)} Given a w.u.C. series $\sum x_n$ in $X$, if the set
$\left\{ P\left( \sum_{i=1}^\infty a_n(i)
x_i\right)\right\}_{n\in\N}$ is relatively weakly compact for
every bounded sequence $(a_n)\subset c_0$, then $\sum x_n$ is u.c.

{\rm (B)} Given a w.u.C. series $\sum x_n$ in $X$, if the set
$\left\{ P\left( \sum_{i=1}^\infty a_n(i)
x_i\right)\right\}_{n\in\N}$ is relatively compact for every
bounded sequence $(a_n)\subset c_0$, then $\sum x_n$ is u.c.

{\rm (C)} If the sequence $(x_n)\subset X$ is equivalent to the
$c_0$-basis, then there is a bounded sequence $(a_n)\subset c_0$
such that the set $\left\{ P\left( \sum_{i=1}^\infty a_n(i)
x_i\right)\right\}_{n\in\N}$ is not relatively weakly compact.

{\rm (D)} If the sequence $(x_n)\subset X$ is equivalent to the
$c_0$-basis, then there is a bounded sequence $(a_n)\subset c_0$
such that the set $\left\{ P\left( \sum_{i=1}^\infty a_n(i)
x_i\right)\right\}_{n\in\N}$ is not relatively compact.

{\rm (E)} For every operator $T\in {\mathcal L}(Z,X)$, if $P\circ
T\in\pkuc Z,Y)$, then $T$ is unconditionally converging.

{\rm (F)} For every operator $T\in {\mathcal L}(c_0,X)$, if
$P\circ T\in\pkK c_0,Y)$, then $T$ is compact.

{\rm (G)} For every subspace $M\subseteq X$ containing a copy of
$c_0$, the polynomial $P\circ j_M$ is not unconditionally
converging, where $j_M$ denotes the embedding of $M$ into $X$.

{\rm (H)} For every subspace $M\subseteq X$ isomorphic to $c_0$,
the polynomial $P\circ j_M$ is not compact, where $j_M$ denotes
the embedding of $M$ into $X$.
\end{Th}

\Proof
(A) $\Ra$ (B) and (C) $\Ra$ (D) are obvious.

(A) $\Ra$ (C): If $(x_n)$ is equivalent to the $c_0$-basis, then
the series $\sum x_n$ is w.u.C., not u.c. So it is enough to apply
(A).

(B) $\Ra$ (D): By the same argument.

(D) $\Ra$ (E): Assume $T\in{\mathcal L}(Z,X)$ is not
unconditionally converging. Then we can find a sequence
$(z_n)\subset Z$ such that $(z_n)$ and $(T(z_n))$ are equivalent
to the $c_0$-basis. By (D), there is a bounded sequence
$(a_n)\subset c_0$ such that the set $\left\{ P\circ T\left(
\sum_{i=1}^\infty a_n(i) z_i\right)\right\}_{n\in\N}$ is not
relatively compact. Hence, letting $M$ be the closed linear span
of $\left\{ z_n\right\}$, there is a bounded sequence
$(x_n)\subset M$ such that $\left\{ P\circ T(x_n)\right\}$ is not
relatively compact. Therefore, $P\circ T\circ j_M$ is not compact.
Since $M$ is isomorphic to $c_0$, this implies that $P\circ T\circ
j_M$ is not unconditionally converging, and we conclude that
$P\circ T$ is not unconditionally converging.

(E) $\Ra$ (F): Take a noncompact operator $T\in{\mathcal
L}(c_0,X)$. Then $T$ is not unconditionally converging. By (E),
$P\circ T$ is not unconditionally converging.

(F) $\Ra$ (G): Suppose there is a subspace $M\subseteq X$
containing $c_0$ such that $P\circ j_M$ is unconditionally
converging. Then, there is a subspace $N\subseteq M$ isomorphic to
$c_0$ so that $P\circ j_N$ is unconditionally converging and so
compact. However, $j_N$ is not compact.

(G) $\Ra$ (H): This is clear, since every unconditionally
converging polynomial on $c_0$ is compact.

(H) $\Ra$ (A): Assume there is a w.u.C. series $\sum x_n$ in $X$,
not u.c., such that the set $\left\{ P\left( \sum_{i=1}^\infty
a_n(i) x_i\right)\right\}_{n\in\N}$ is relatively weakly compact
for every bounded sequence $(a_n)\subset c_0$. Taking blocks, we
can assume that $(x_n)$ is equivalent to the $c_0$-basis. Let $M$
be the closed linear span of $\left\{ x_n\right\}$. Then, $P\circ
j_M$ takes bounded sequences into relatively weakly compact
sequences, and so $P\circ j_M$ is compact.\finesp

\begin{Def}{\rm
We say that $P\in\pkXY$ is an {\it Orlicz-Pettis polynomial\/} if
it satisfies the equivalent  assertions of Theorem~\ref{uc+}. We
denote by $\pkucmXY$ the space of all $k$-homogeneous
Orlicz-Pettis polynomials from $X$ into $Y$.}
\end{Def}

The choice of the name is due to the relationship with the u.c.\
series, which were studied by Orlicz and Pettis
\cite[Chapter~IV]{Di}.

The classes $\pkucXY$ and $\pkucmXY$ may be described by means of
a family of sets. We say that a subset $A\subset X$ is a {\it
WUC-set\/} if there is an operator $T\in{\mathcal L} (c_0,X)$ such
that $A=T\left(B_{c_0}\right)$. We need a previous lemma.

\begin{Lemma}
Given a polynomial $P\in\pkXY$ which is not unconditionally
converging, there is an embedding $j:c_0\ra X$ such that $P\circ
j\notin\pkuc c_0,Y)$.
\end{Lemma}

\Proof
If $P\notin \pkucXY$, we can find a w.u.C. series $\sum x_n$ in
$X$ such that the sequence $\left( P(x_1+\cdots +x_n)\right)_n$ is
not convergent. Let $j:c_0\ra X$ be given by $j(e_n)=x_n$. Then
the sequence
$$
\left( P\circ j(e_1+\cdots +e_n)\right)_n=
\left( P(x_1+\cdots +x_n)\right)_n
$$
is not convergent, so $P\circ j\notin\pkuc c_0,Y)$.\finesp

The next Proposition highlights the opposition between the classes
$\pkucXY$ and $\pkucmXY$.

\begin{Prop}
For a polynomial $P\in\pkXY$, we have:

{\rm (a)} $P\in\pkucXY$ if and only if $P$ takes WUC-sets into
relatively (weakly) compact sets;

{\rm (b)} $P\in\pkucmXY$ if and only if, for every WUC-set
$A\subset X$, if $P(A)$ is relatively (weakly) compact, then so is
$A$.
\end{Prop}

\Proof
(a) Let $P\in\pkucXY$ and $A=T\left(B_{c_0}\right)$. Then $P\circ
T\in\pkuc c_0,Y)$, and so $P(A)=P\circ T\left(B_{c_0}\right)$ is
relatively compact. Conversely, suppose $P\notin\pkucXY$. By the
Lemma, there is $T\in{\mathcal L}(c_0,X)$ such that $P\circ
T\notin\pkuc c_0,Y)$. Therefore, $T\left(B_{c_0}\right)$ is a
WUC-set so that $P\circ T\left(B_{c_0}\right)$ is not relatively
weakly compact.

(b) Let $P\in\pkucmXY$ and choose a WUC-set $A\subset X$ so that
$P(A)$ is relatively weakly compact. Take $T\in{\mathcal
L}(c_0,X)$ with $A=T\left(B_{c_0}\right)$. Then $P\circ
T\left(B_{c_0}\right)$ is relatively weakly compact, so $P\circ T$
is compact. Since $P\in\pkucmXY$, $T$ is compact and $A$ is
relatively compact. Conversely, let $T\in{\mathcal L}(c_0,X)$ with
$P\circ T$ compact. Since $P\circ T\left(B_{c_0}\right)$ is
relatively compact, we have that $T\left(B_{c_0}\right)$ is
relatively compact, so $T$ is compact. Therefore,
$P\in\pkucmXY$.\finesp

The following result gives a polynomial satisfying the assertions
of Theorem~\ref{uc+}. Other examples are shown in
Section~\ref{counterex}.

\begin{Prop}
\label{gammak}
For every Banach space $X$, the polynomial $\gamma_k:X\ra\ptpp^k
X$ given by
$$
\gamma_k(x)= x^{(k)}
$$
is an Orlicz-Pettis polynomial.
\end{Prop}

\Proof
Take a subspace $M\subseteq X$ isomorphic to $c_0$, and let
$j:c_0\ra M$ be a surjective isomorphism. From
$$
\left\|\gamma_k\circ j_M\circ j(e_n)\right\| =
\left\| j_M\circ j(e_n)\right\|^k\not\longrightarrow 0,
$$
we get that $\gamma_k\circ j_M$ is not compact. Then apply (H) of
Theorem~\ref{uc+}.\fin

\section{Counterexamples}
\label{counterex}

In this Section, we first give some more properties of the
polynomial $\gamma_k$ considered in Proposition~\ref{gammak}, some
of which are used to establish a theorem about polynomials on
spaces containing $c_0$. All these previous results are applied in
the main theorem of the Section that provides sufficient
conditions for a polynomial to be Orlicz-Pettis. A number of
counterexamples are given to show that these conditions are not
necessary.

Our first theorem gives a property of polynomials on spaces
containing a copy of $c_0$. We need two previous results.

\begin{Lemma}
\label{c0basis}
Given $k\in \N$, the sequence $\left( e_n^{(k)}\right)_n$ in
$\ptpp^k c_0$ is equivalent to the unit vector basis in $c_0$.
\end{Lemma}

\Proof
By induction on $k$, we show that
$$
\left\| a_1e_1^{(k)}+\cdots+a_ne_n^{(k)}\right\| =\left\|
a_1e_1+\cdots+a_ne_n\right\|
$$
where $(a_i)_{i=1}^n$ is a finite sequence of scalars.

For $k=1$ there is nothing to prove. Suppose the result holds for
$k-1$, and let $r_n(t)=\sign\sin 2^n\pi t$ for $t\in [0,1]$. Then,
assuming $\max |a_i|=1$, we have
\be
\lefteqn{ a_1e_1^{(k)}+\cdots+a_ne_n^{(k)}=}\\
&=&\int_0^1\left[ a_1r_1(t)e_1^{(k-1)}+\cdots+a_nr_n(t)e_n^{(k-1)}
    \right]\otimes\left[ a_1r_1(t)e_1+\cdots+a_nr_n(t)e_n\right]
    dt\\
&=&2^{-n}\sum_{i=1}^{2^n}\left[ a_1\eps_1(i)e_1^{(k-1)}+\cdots+
    a_n\eps_n(i)e_n^{(k-1)}\right]\otimes\left[ a_1\eps_1(i)e_1+
    \cdots+a_n\eps_n(i)e_n\right]
\ee
where $\eps_j(i)$ is the value of $r_j(t)$ on the interval
$$
\left(\frac{i-1}{2^n}\; ,\,\frac{i}{2^n}\right)\qquad\mbox{for
$1\leq i\leq 2^n$} .
$$
By the induction hypothesis, we get
$$
\left\| a_1e_1^{(k)}+\cdots+a_ne_n^{(k)}\right\|\leq 1.
$$

On the other hand, there is $i_0\in\{ 1,\ldots,k\}$ so that
$$
\left\| a_1e_1+\cdots+a_ne_n\right\|=\left| a_{i_0}\right| .
$$
Considering $e_{i_0}$ as a vector of $\ell_1$, take
$$
e_{i_0}^{(k)}\in\left( \ptpp^kc_0\right) ^*\approx\Lk c_0).
$$
Clearly,
$$
\left\| e_{i_0}^{(k)}\right\|=1,\quad\mbox{and}\quad
\left\langle e_{i_0}^{(k)},a_1e_1^{(k)}+\cdots+a_ne_n^{(k)}
\right\rangle = a_{i_0},
$$
so the result follows.\fin

\begin{Prop}
\label{wuCucc0}
There is a w.u.C. series $\sum y_i$ in $c_0$, not u.c., such that
$\sum_i\gamma_k(y_i)$ is u.c. in $\ptpp^kc_0$ for each $k\geq 2$.
\end{Prop}

\Proof
For simplicity, consider the case $k=2$. Take the vectors
$$
y_i=\frac{e_n}{n}\quad\mbox{for}\quad n\in\N\; ,\quad
\frac{n(n-1)}{2}<i\leq\frac{n(n+1)}{2}\; .
$$
Clearly, the series $\sum y_i$ is w.u.C\@. Since
$$
\left\|\sum_{i=1+n(n+1)/2}^{m(m+1)/2} y_i\right\| =1\quad
\mbox{for}\quad m>n,
$$
the series is not u.c. Moreover, for every finite sequence
$$
\frac{n(n+1)}{2}<i_1<\cdots <i_l,
$$
we have, from Lemma~\ref{c0basis},
$$
\left\|\sum_{j=1}^{l} y_{i_j}\otimes y_{i_j}\right\| <\frac1n\; .
$$
Therefore, $\sum y_i\otimes y_i$ is u.c. in $c_0\ptpp c_0$.\fin

\begin{Th}
\label{copc0tri}
Given Banach spaces $X$ and $Y$, with $X$ containing a copy of
$c_0$, an integer $k>1$ and a polynomial $P\in\pkXY$, we can find
a w.u.C. series $\sum x_i$ in $X$, not u.c., such that $\sum
P(x_i)$ is u.c. in $Y$.
\end{Th}

\Proof
Let $j:c_0\ra X$ be an embedding. Consider the commutative diagram
$$
\begin{CD}
 X         @>\gamma_k>>       \ptpp^k X\\
 @AjAA                        @AA\otimes^kjA\\
 c_0       @>>\gamma_k>       \ptpp^kc_0
\end{CD}
$$

Let $\sum y_i$ be the series constructed in the proof of
Proposition~\ref{wuCucc0}. Then, the series $\sum j(y_i)$ is
w.u.C., not u.c., in $X$, and the series $\sum_i\gamma_k\circ
j(y_i)=\sum_i\left(\otimes ^k j\right)\circ \gamma_k(y_i)$ is u.c.
Let $\tilde{P}:\ptpp^kX\ra Y$ be the operator defined by
$\tilde{P}(x_1\otimes\cdots\otimes x_k):=\hat{P}(x_1,\ldots,x_k)$,
where $\hat{P}\in\LkXY$ is the symmetric $k$-linear mapping
associated to $P$. Then the series $\sum_i P\circ
j(y_i)=\sum_i\tilde{P}\circ
\gamma_k\circ j(y_i)$ is u.c. in $Y$.\finesp

The next theorem shows that the polynomial $\gamma_k:X\ra\ptpp^kX$
takes sequences equivalent to the $c_0$-basis into sequences
equivalent to the $c_0$-basis. Again, we need two preparatory
results.

\begin{Lemma}
\label{gammac0}
The polynomial $\gamma_k:c_0\ra\ptpp^k c_0$ takes sequences
equivalent to the $c_0$-basis into sequences equivalent to the
$c_0$-basis.
\end{Lemma}

\Proof
Letting $j:c_0\ra c_0$ be an isomorphism, consider the commutative
diagram
$$
\begin{CD}
 c_0       @>\gamma_k>>       \ptpp^k c_0\\
 @AjAA                        @AA\otimes^kjA\\
 c_0       @>>\gamma_k>       \ptpp^kc_0
\end{CD}
$$

Since $\left(\gamma_k(e_n)\right)_n$ is equivalent to the
$c_0$-basis (Lemma~\ref{c0basis}), it is enough to show that
$\otimes^kj$ is an injective isomorphism. Since $j(c_0)$ is
complemented in $c_0$, there is an operator $S:c_0\ra c_0$ such
that $S\circ j=I_{c_0}$. Then, $\left(\otimes^kS\right)\circ
\left(\otimes^kj\right)$ is the identity map on $\ptpp^kc_0$.
Hence, $\otimes^kj$ is an injective isomorphism.\fin

\begin{Prop}
Let $j:c_0\ra X$ be an injective isomorphism. Then the operator
$$
j^{(k)}:\ptpp^kc_0\lra\ptpp^kX
$$
is an injective isomorphism.
\end{Prop}

\Proof
Take $z\in\ptpp^kc_0$ with $z=\sum_{i=1}^\infty x^1_i\otimes\cdots
\otimes x^k_i$, and $A\in\left(\ptpp^kc_0\right)^*\approx\Lk c_0)$,
with $\| A\|=1$ and $\langle A,z\rangle=\|z\|$. There is an
extension $\tilde{A}\in\Lk\ell_\infty)$ of $A$ with
$\|\tilde{A}\|=\|A\|$
\cite{GV}. Consider the second adjoint $j^{**}:\ell_\infty\ra X^{**}$
of $j$. By the injectivity of $\ell_\infty$, the operator
$(j^{**})^{-1}:j^{**}(\ell_\infty)\ra\ell_\infty$ has an extension
to an operator $\pi:X^{**}\ra\ell_\infty$; clearly, $\pi\circ
j^{**}=I_{\ell_\infty}$. Let $B:=\tilde{A}\circ (\pi\circ
J_X)^k\in\LkX$, where $J_X:X\ra X^{**}$ is the canonical
embedding. Then, $\|B\|\leq\|\tilde{A}\|\cdot\|\pi\|^k$ and
$$
B\circ j^k=\tilde{A}\circ (\pi\circ J_X\circ j)^k=
\tilde{A}\circ \left(\pi\circ j^{**}\circ J_{c_0}\right)^k=
\tilde{A}\circ J_{c_0}^k=A.
$$
Therefore,
\be
\|\pi\|^k\cdot\|A\|\cdot\|j^{(k)}(z)\|&\geq&|\langle
B,j^{(k)}(z)\rangle|\\
&=&\left|\sum_{i=1}^\infty B\left( j\left(x^1_i\right) ,\ldots,
j\left(x^k_i\right)\right)\right|\\
&=&\left|\sum_{i=1}^\infty A\left( x^1_i,\ldots,x^k_i\right)
\right|\\
&=&|\langle A,z\rangle|\\
&=&\|z\| ,
\ee
and this finishes the proof.\fin

\begin{Th}
The polynomial $\gamma_k:X\ra\ptpp^kX$ takes sequences equivalent
to the $c_0$-basis into sequences equivalent to the $c_0$-basis.
\end{Th}

\Proof
If $X$ contains no copy of $c_0$, the result is trivially true. If
$X=c_0$, see Lemma~\ref{gammac0}. If $X$ contains a copy of $c_0$,
the last Proposition reduces the problem to the case
$X=c_0$.\finesp

The following result gives an example of a polynomial on $c_0$
which will be useful.

\begin{Prop}
\label{Pen0}
There is a polynomial $P\in{\mathcal P}(^2c_0,c_0)$ such that
$P(e_n)=0$ for all $n$, but $P$ is not compact on any infinite
dimensional subspace.
\end{Prop}

\Proof
Consider a bijection
$$
(\alpha,\beta):\N\lra\{ (n,m)\in\N\times\N:n\neq m\} .
$$
Define
$$
P(x):=\left( x(\alpha (i))x(\beta (i))\right)_{i=1}^\infty\quad
\mbox{for}\quad x=(x(i))\in c_0.
$$
Then $P(e_n)=0$ for all $n$. If $M\subseteq c_0$ is an infinite
dimensional subspace, we can find a norm one sequence
$(x_n)\subset c_0$ disjointly supported such that $\dist
(x_n,M)<2^{-n}$. For each $n\in\N$, let $k_n\in\N$ satisfy
$|x_n(k_n)|=1$. Then
$$
\|P(x_{2n}+x_{2n+1})\|\geq |x_{2n}(k_{2n})x_{2n+1}(k_{2n+1})|=1,
$$
which implies that $P$ is not compact on  $M$.\finesp

We can now state the main result of the Section.

\begin{Th}
\label{contra}
Let $P\in\pkXY$ be a polynomial, with $k\geq 2$. Consider the
following assertions:

{\rm (A)} If $\sum x_n$ is w.u.C. in $X$, and $\sum P(x_n)$ is
u.c. in $Y$, then $\sum x_n$ is u.c.

{\rm (B)} Every sequence $(x_n)\subset X$ equivalent to the
$c_0$-basis has a subsequence $\left( x_{n_i}\right)$ such that
$\left( P\left( x_{n_i}\right)\right)$ is equivalent to the
$c_0$-basis.

{\rm (C)} If $\sum x_n$ is w.u.C. in $X$, and $\left(
P\left(\sum_{i=1}^n x_i\right)\right)_n$ is convergent, then $\sum
x_n$ is u.c.

{\rm (D)} If the sequence $(x_n)\subset X$ is equivalent to the
$c_0$-basis, then $\left( P\left(\sum_{i=1}^n x_i\right)\right)_n$
is not relatively compact.

{\rm (E)} If the sequence $(x_n)\subset X$ is equivalent to the
$c_0$-basis, then $\lim \|P(x_n)\|\not\ra 0$.

{\rm (F)} $P$ is an Orlicz-Pettis polynomial.

Then the following and only the following implications hold: {\rm
$$
\begin{array}{ccccc}
 \mbox{(A)} & \Longrightarrow & \mbox{(B)} & \Longleftrightarrow &
 \mbox{(E)} \\
         \Downarrow   &  &            &  &      \Downarrow      \\
 \mbox{(C)} & \Longleftrightarrow & \mbox{(D)} & \Longrightarrow &
 \mbox{(F)} \\
\end{array}
$$}
\end{Th}

\Proof
(A) $\Ra$ (B) and (A) $\Ra$ (C): If $P$ satisfies (A), then
Theorem~\ref{copc0tri} implies that $X$ contains no copy of $c_0$.
So (B) and (C) are satisfied in a trivial way.

(B) $\Ra$ (E) is obvious.

(E) $\Ra$ (B): Let $(x_n)\subset X$ be equivalent to the
$c_0$-basis. Then $\sum x_n$ is w.u.C., so $\sum P(x_n)$ is also
w.u.C. \cite{GG}. In particular, $(P(x_n))$ is weakly null. By
(E), passing to a subsequence, we can assume that $(P(x_n))$ is
seminormalized and basic, so it is equivalent to the $c_0$-basis
\cite[Corollary~V.7]{Di}.

(B) $\not\Ra$ (A): Consider the polynomial
$\gamma_k:c_0\ra\ptpp^kc_0$ (see Proposition~\ref{wuCucc0} and
Lemma~\ref{gammac0}).

(D) $\not\Ra$ (A): Consider the polynomial defined in
Proposition~\ref{Pen0}.

(C) $\Ra$ (D): Assume $P$ does not satisfy (D). Then there is a
sequence $(x_n)\subset X$ equivalent to the $c_0$-basis, such that
$\left( P\left(\sum_{i=1}^n x_i\right)\right)_n$ is relatively
compact. We can find an increasing sequence of indices $(m_i)$ so
that $\left( P\left(\sum_{i=1}^n y_i\right)\right)_n$ is
convergent, where
\begin{equation}
\label{block}
y_i=\sum_{j=m_i+1}^{m_{i+1}} x_j\, .
\end{equation}
Since $(y_i)$ is equivalent to the $c_0$-basis, $P$ does not
satisfy (C).

(D) $\Ra$ (C): Assume $P$ does not satisfy (C). Then there is a
w.u.C. series $\sum x_n$ in $X$, not u.c., so that $\left(
P\left(\sum_{i=1}^n x_i\right)\right)_n$ is convergent. Take an
increasing sequence of indices $(m_i)$ so that $(y_i)$ is
equivalent to the $c_0$-basis, where $y_i$ is defined as in
(\ref{block}). Then $\left( P\left(\sum_{i=1}^n
y_i\right)\right)_n$ is a subsequence of $\left(
P\left(\sum_{i=1}^n x_i\right)\right)_n$ and so it converges, in
contradiction with (D).

(E) $\Ra$ (F): Assume $T\in{\mathcal L}(Z,X)$ is not
unconditionally converging. Then we can find a sequence
$(z_n)\subset Z$ such that $(z_n)$ and $(T(z_n))$ are equivalent
to the $c_0$-basis. If $P$ satisfies (E), we have $\|P\circ
T(z_n)\|\not\to 0$, which implies that $P\circ T\not\in\pkuc
Z,Y)$. So, by Theorem~\ref{uc+}(E), $P\in\pkucmXY$.

(F) $\not\Ra$ (E): The polynomial $P$ of Proposition~\ref{Pen0}
does not satisfy (E). To see that it does satisfy (F), take an
operator $T\in{\mathcal L}(Z,c_0)$ not unconditionally converging.
There is an operator $j:c_0\ra Z$ such that $(T\circ j(e_n))$ is
equivalent to $(e_n)$. Passing to a perturbed subsequence, we can
assume that $(T\circ j(e_n))$ is disjointly supported. The series
$\sum_n\left( j(e_{2n})+j(e_{2n+1})\right)$ is w.u.C\@. However,
$\left\| P\circ T\left(j(e_{2n})+j(e_{2n+1})\right)\right\|$ is
bounded away from $0$, so $P\circ T$ is not unconditionally
converging. By Theorem~\ref{uc+}(E), $P\in\pkucm Z,c_0)$.

(D) $\Ra$ (F): Assume $T\in{\mathcal L}(Z,X)$ is not
unconditionally converging. Then we can find a sequence
$(z_n)\subset Z$ such that $(z_n)$ and $(T(z_n))$ are equivalent
to the $c_0$-basis. If $P$ satisfies (D), the sequence $\left(
P\circ T\left(\sum_{i=1}^n z_i\right)\right)_n$ is not relatively
compact. Hence, $P\circ T\not\in \pkucXY$. By
Theorem~\ref{uc+}(E), $P\in\pkucmXY$.

(D) $\not\Ra$ (E): Let $P$ be the polynomial defined in
Proposition~\ref{Pen0}, and $(x_n)\subset c_0$ a sequence
equivalent to the $c_0$-basis. Denote $y_n:=x_1+\cdots +x_n$, and
$z:=\sum_{n=1}^\infty x_n\in\ell_\infty\backslash c_0$. Let
$$
3\delta :=\limsup_i z(i)>0.
$$
If the $\limsup$ were not positive, then we would take the
$\liminf$. Choose $i_1\in\N$ with $|z(i_1)-3\delta|<\delta$. There
is $n_1\in\N$ so that $|y_n(i_1)-3\delta|<\delta$ for all $n\geq
n_1$. Choose now $i_2\in\N$ $(i_2>i_1)$ so that $\left|
y_{n_1}(i)\right|<\delta/2$ for all $i\geq i_2$ and
$|z(i_2)-3\delta|<\delta$. There is $n_2\in\N$ $(n_2>n_1)$ so that
$\left| y_n(i_2)-3\delta\right|<\delta$ for all $n\geq n_2$.

Proceeding in this way, we obtain two increasing sequences of
integers $(i_j)$, $(n_j)$ so that, for $j<l$,
\be
\left\| P\left( y_{n_j}\right)-P\left( y_{n_l}\right)\right\|
&\geq& \left| y_{n_j}(i_1)y_{n_j}(i_l)-y_{n_l}(i_1)y_{n_l}(i_l)
\right|\\
&\geq&\left| y_{n_l}(i_1)\right|\cdot\left| y_{n_l}(i_l)\right| -
\left| y_{n_j}(i_1)\right|\cdot\left| y_{n_j}(i_l)\right|\\
&>&2\delta^2.
\ee
Therefore, $P$ satisfies (D). Clearly, $P$ does not satisfy (E).

(E) $\not\Ra$ (D): Let $P\in{\mathcal P}(^2\!c_0,c_0)$ be given by
$$
P(x)=\sum_{j=2}^\infty\sum_{i=1}^{j-1}\left( x(j)-x(i)\right)
x(j)e_{j^2+i}\qquad\mbox{for $x=(x(i))_{i=1}^\infty\in c_0$}.
$$
Clearly, $P$ satisfies (E). Since $P(e_1+\cdots +e_n)=0$ for all
$n$, $P$ does not satisfy (D).\fin

\begin{Rem}{\rm
In the linear case $(k=1)$, all the assertions of
Theorem~\ref{contra} are equivalent
\cite{GM}. So, our choice for the definition of the Orlicz-Pettis
polynomials provides the widest possible class.}
\end{Rem}

\end{document}